\documentclass [12pt,a4paper]{article}
\usepackage{epsfig,amsmath,amssymb,color}
\usepackage{natbib}
\usepackage{bm,bbm}
\usepackage{psfrag}
\usepackage{wrapfig}
\usepackage{geometry}
\newcommand{\bx}{\mathbf{x}}
\newcommand{\bX}{\mathbf{X}}

\newtheorem{thm}{Theorem}[section]
\newtheorem{definition}[thm]{Definition}
\newtheorem{proposition}[thm]{Proposition}

\def\e{{\mathbb{E}}}
\def\prob{{\mathbb{P}}}

\def\convp{\stackrel{{\scriptsize{\prob}}}{\longrightarrow}}
\def\convlr{\stackrel{\mbox{\scriptsize{$L_r$}}}{\longrightarrow}}
\def\convltwo{\stackrel{\mbox{\scriptsize{$L_2$}}}{\longrightarrow}}

\def\qed{\hfill\hbox{${\vcenter{\vbox{
        \hrule height 0.4pt\hbox{\vrule width 0.4pt height 6pt
        \kern5pt\vrule width 0.4pt}\hrule height 0.4pt}}}$}}

\def\Ref#1{(\ref{#1})}

\def\var{\hbox{Var}}

\def\min{\hbox{min}}

\def\cdf{\hbox{cdf\,}}

\def\iid{\hbox{i.i.d.}}
\def\pdQ{\hbox{pdQ}}
\def\pdQs{\hbox{pdQs}}

\newcommand{\F}{\mathcal{F}\,}

\begin{document}

\title{Divergence from, and Convergence to, Uniformity of Probability Density Quantiles}

\author{Robert G. Staudte\footnote{ Postal address: Department of Mathematics and Statistics, La Trobe University, VIC 3086, Australia. Email address: r.staudte@latrobe.edu.au.}\\ {\it La Trobe University}
\\ Aihua Xia\footnote{Postal address:
School of Mathematics and Statistics, University of Melbourne,
VIC 3010, Australia. Email address: aihuaxia@unimelb.edu.au. Research supported by ARC Discovery Grant DP150101459.} \\ {\it University of Melbourne}
}

\date{15 March, 2018}
\maketitle

\abstract{We demonstrate that questions of convergence, divergence and inference regarding shapes of distributions
can be carried out in a location- and scale-free environment. This environment is the class of probability density quantiles (\pdQs), obtained by normalizing the composition of the density with the associated quantile function. It has earlier been shown that the \pdQ\  is representative of a  location-scale family and carries essential information regarding shape and tail behavior of the family. The class of \pdQs \  are densities of continuous distributions with common domain, the unit interval, facilitating
 metric and semi-metric comparisons. Further applications of the \pdQ\ mapping are quite generally entropy increasing so convergence to the uniform distribution is investigated. New fixed point theorems are established and illustrated by examples.   The Kullback-Leibler directed divergences from uniformity of these \pdQs\ are mapped and  found to be essential ingredients in power functions of optimal tests for uniformity against alternative shapes.}

{\em Keywords:\ convergence in {$L_r$ norm}; fixed point theorem;  Kullback-Leibler divergence; relative entropy;
semi-metric; uniformity test}

\section{Introduction}\label{sec:intro}

\subsection{Background and summary}\label{sec:background}

For each continuous location-scale family of distributions with square-integrable density there is a probability density quantile (\pdQ) which is an absolutely continuous distribution on the unit interval. Members of the class of such \pdQs\ differ only in {\em shape}, and the asymmetry of their  shapes can be partially ordered by
their Hellinger distances or Kullback-Leibler divergences from the class of symmetric distributions on this interval.
In addition, the tail behaviour of the original family can be described in terms of the boundary derivatives of its \pdQ. Empirical estimators of the \pdQs\ enable one to carry out inference, such as fitting
shape parameter families to data; details are in \cite{S-2017}.

The Kullback-Leibler directed divergence and symmetrized divergence (KLD) of a \pdQ\ with respect to the uniform distribution on [0,1] is investigated in Section~\ref{sec:defns}, with remarkably simple numerical results, and
a map of these divergences for some standard location-scale families is constructed. The \lq shapeless\rq\ uniform distribution is the center of the \pdQ\ universe, as is explained in Section~\ref{sec:convto}, where it is found to be a fixed point. We then investigate the convergence to uniformity of repeated applications of the \pdQ\ transformation, by means of
fixed point theorems for a semi-metric. In Section~\ref{sec:testing} power functions of hypothesis tests of uniformity against specified alternative shapes are shown to be dependent on the symmetrized Kullback-Leibler divergence.
 Further ideas are discussed in Section~\ref{sec:summary}.

\subsection{Definitions and divergence map}\label{sec:defns}

 Let $\F $ denote the class of {cumulative distribution functions (\cdf s)  on the real line} and for each $F\in \F$ define the associated {\em quantile function} of $F$ by $ Q(u)=\inf \{x:\ F(x)\ge u \}$, for $0< u< 1.$ When the random variable $X$ has \cdf \,$F$, we write $X\sim F$. {When the density function $f=F'$ exists, we also write
 $X\sim f$ or $f\sim F$.} We only discuss $F$ absolutely continuous with respect to Lebesgue measure, but the results can be extended to the discrete and mixture cases using {suitable} dominating measures.

\begin{definition}\label{def1}
Let $\F '=\{F\in \F :\; f=F' \text { exists and is positive}\}$. For each $F\in \F '$ we follow \cite{par-1979} and define  the {\em quantile density} function $q(u)=Q'(u)=1/f(Q(u))$. Parzen called its reciprocal function $fQ(u)= f(Q(u))$  the {\em density quantile} function.
  For $F\in \F '$, and $U$ uniformly distributed on [0,1], assume $\kappa =\e [fQ(U)] = \int f^2(x)\,dx $
 is finite; that is, $f$ is square integrable.
 Then we can define the {\em continuous {\pdQ }} of $F$ by  $f^*(u)=fQ(u)/\kappa $, $0<u<1$.  Let $\F '^{*}\subset \F ' $ denote the class of all such $F$.
\end{definition}

Not all $f$  are square-integrable, and this requirement for the mapping $f \to f^*$ means that
$\F '^{*}$ is a proper subset of $\F ' .$ The advantages of working with $f^*$s over $f$s are that they are free of location and scale parameters, they ignore flat spots in $F$ and have a common bounded support. Moreover, $f^*$ often has a simpler formula than $f$; see Table~\ref{table1} for examples.

Next we evaluate and plot the \cite{kl-1951}  divergences from uniformity.
The \cite{kl-1951} divergence of density $f_1$ from density $f_2$, {when both have domain [0,1], is defined as
$I(f_1:f_2) := \int _0^1\ln (f_1(u)/f_2(u))\,f_1(u)\,du =\e[\ln (f_1(U)/f_2(U))\,f_1(U)],$} where $U$ denotes a random variable
with the uniform distribution ${\cal U}$ on [0,1].
The divergences from uniformity are easily computed through
$I({\cal U}:f^*) = -\int _0^1\ln (f^*(u))\,du =-\e[\ln (f^*(U))]$ and
$I(f^*:{\cal U}) = \int _0^1\ln (f^*(u))\,f^*(u)\,du = \e[\ln (f^*(U))\,f^*(U)].$
 \cite[p.~6]{Kullback-1968} interprets $I(f^*:{\cal U})$ {as} the mean evidence in one observation $V\sim f^*$ for
$f^*$ over ${\cal U}$; it is also known as the {\em relative entropy} of $f^*$ with respect to
$ {\cal U}$.  In Table~1 are shown the quantile functions of some standard  distributions, along with their \pdQs,\ associated divergences $I({\cal U}:f^*),I(f^*:{\cal U})$ and
    symmetrized divergence (KLD)\  defined by $J({\cal U},f^*) := I({\cal U}:f^*)+I(f^*:{\cal U})  $.

 \begin{figure}[t!]
\begin{center}
\includegraphics[scale=.8]{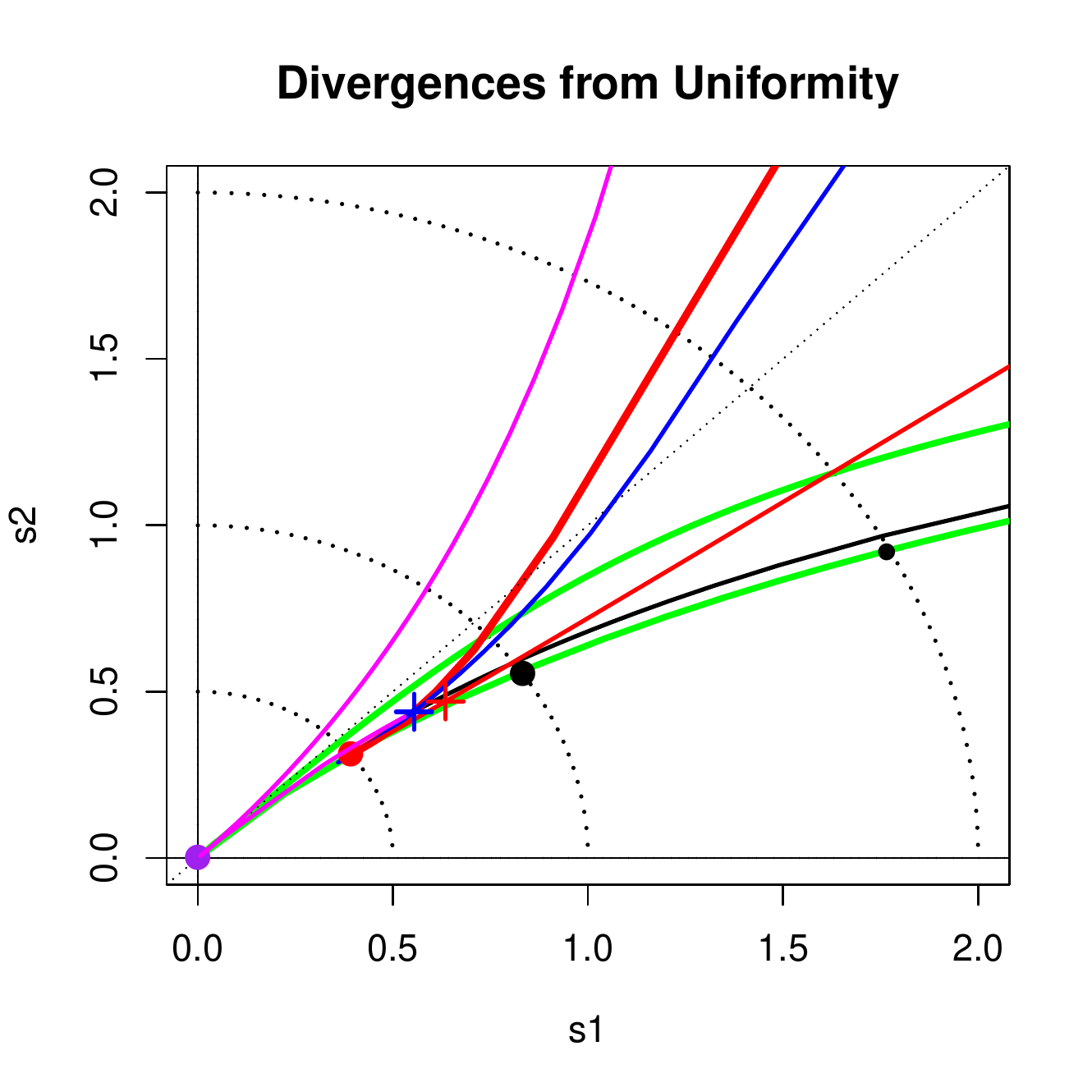}
\caption{\footnotesize  {\bf Divergence from uniformity.} \em The loci of points $(s_1,s_2)=(\sqrt {I({\cal U}:f^*)}\,,\sqrt {I(f^*:{\cal U})}\,)$ is shown for various standard families.
The large disks correspond respectively to the symmetric families:\; uniform (purple), normal (red) and Cauchy (black). The crosses correspond to the asymmetric distributions:\; exponential (blue) and standard lognormal (red).   More details are given in Section~\ref{sec:defns}.\label{fig1}}
\end{center}
\end{figure}

\begin{definition}\label{def2}
Given \pdQs\  $f_1^*$, $f_2^*$, let $d(f_1^*,f_2^*):=\sqrt {I(f^*_1:f^*_2)+I(f_2^*:f_1^*) }$. {Then $d$ is a semi-metric on the space of \pdQs; i.e., $d$ satisfies all requirements of a metric except the triangle inequality. Introducing  the  coordinates $(s_1,s_2)=(\sqrt {I({\cal U}:f^*)}\,,\sqrt {I(f^*:{\cal U})})$, we can define the
{\em distance from uniformity} of any $f^*$ by the  Euclidean distance of $(s_1,s_2)$ from the origin $(0,0)$, namely $d({\cal U},f^*)$.}
\end{definition}

{{\bf Remark} This $d$ does not satisfy the triangle inequality:\
for example, if $\cal U,$ $\cal N$ and $\cal C$ denote the uniform, normal and Cauchy \pdQs, then
$d({\cal U},{\cal N})=0.5,$  $ d({\cal N},{\cal C})=0.4681$ but $d({\cal U},{\cal C})=1.$ However, $d$
can provide an informative measure of distance from uniformity.}

In Figure~\ref{fig1} are shown the loci of points $(s_1,s_2)$ for some continuous shape families. The light dotted arcs with radii 1/2, 1 and 2 are a guide to these distances from uniformity.
The large discs in purple, red and black correspond to  $\cal U,$ $\cal N$ and $\cal C$.
The blue cross at distance $1/\sqrt {2}$ from the origin corresponds to the exponential distribution.  Nearby is the standard lognormal point marked by a red cross. The lower red curve is nearly straight and is
the locus of points corresponding to the lognormal shape family.

 The Chi-squared($\nu$), $\nu >1$, family also appears as a red curve; it passes through the blue cross when $\nu =2$,
 as expected, and heads toward the normal disc as $\nu \to \infty .$  The Gamma family has the same locus of points
 as the Chi-squared family.
  The curve for the Weibull($\beta $) family,  for $0.5 < \beta <3${,} is shown in blue; it crosses the exponential blue cross when $\beta =1$. The Pareto($a$)  curve is shown in black. As $a$ increases from 0, this line crosses the arcs  distant 2 and 1 from the origin for $a=(2\sqrt{2}\,+1)/7\approx 0.547$ and $a=(\sqrt {5}\,-1)/2\approx 1.618$, respectively, and approaches the exponential blue cross as $a\to \infty $.

The Power($b$) or Beta($b,1$) for $b>1/2$ family is represented by the magenta curve of points moving toward the origin as $b$ increases from 1/2 to 1, and then moving
out towards the exponential blue cross as $b\to \infty $. For each choice of $\alpha >0.5,$ $\beta >0.5$ the locus of the Beta($\alpha,\beta $)
\pdQ\ divergences lies above the chi-squared  red curve and mostly below the power($b$) magenta curve; however, the U-shaped Beta distributions have loci above it.

The lower green line near the Pareto black curve gives the loci of root-divergences from uniformity of the Tukey($\lambda $) with $\lambda <1$, while the upper green curve
corresponds to $\lambda \geq 1$.
 It is known that the Tukey($\lambda$) distributions, with $\lambda <1/7$, are good approximations to  Student's $t$ distributions for $\nu >0$
provided $\lambda $ is chosen properly. The same is true for their corresponding \pdQs\ \cite[Sec.3.2]{S-2017}.
For example, the \pdQ\ of $t_\nu $ with $\nu =0.24$ degrees of freedom is well approximated by the choice $\lambda =-4.063.$ Its location is marked by the small black disk
in Figure~\ref{fig1}; it is distant 2 from uniformity.
The generalized Tukey distributions of \cite{fmkl-1988}  with two shape parameters also fill a large funnel
shaped region (not marked on the map) emanating from the origin and just including the region bounded by the
green curves of the Tukey symmetric distributions.

The larger the value of $d({\cal U},f^*)$, the easier it should be to discriminate between $\cal U$ and $f^*$
based on data, as we will see when calculating power functions of tests in Section~\ref{sec:power}.
First, however, we show that $\cal U$ plays a unique role in the space of \pdQs.

\section{{Convergence of density shapes to uniformity via fixed point theorems}}\label{sec:convto}

The transformation $f\to f^*$ of Definition~\ref{def1} is quite powerful, removing location and scale
and moving the distribution from the support of $f$ to the unit interval. Examples suggest
 that another application of the transformation $f^{2*}:=(f^*)^*$ leaves less information about $f$ in $f^{2*}$ and hence it is closer to the uniform density.  Further, with $n$ iterations $f^{(n+1)*}:=(f^{n*})^*$ for $n\geq 2$, we would expect that $f^{n*}$ converges to the uniform density as $n\to \infty$. An R script \cite{R} for finding repeated $*$-iterates of a given \pdQ\; is available as Supplementary Online Material.

\subsection{Conditions for convergence to uniformity}\label{sec:convunif}

\begin{definition}\label{def3}
Given $f\in \F '$,  we say that $f$ is of {\em $*$-order $n$} if $f^{*},f^{2*},\dots ,f^{n*}$ exist
but $f^{(n+1)*}$ does not.  When the infinite sequence $\{f^{n*}\}_{n\geq 1}$ exists, it is said to be of {\em infinite $*$-order}.
\end{definition}
For example, the Power($3/4$) family is of $*$-order $2$, while the  Power($2$) family is of infinite
$*$-order. The $\chi ^2_\nu $ distribution is of finite $*$-order  for $1<\nu <2$ and infinite $*$-order for $\nu \geq 2.$ The normal distribution is of infinite $*$-order.

\bigskip
We write $\mu_n:=\int_{-\infty}^\infty \{f(y)\}^n\,dy$, $\kappa_n=\int_0^1 \{f^{n*}(x)\}^2\,dx$, $n\ge 1$, and $\kappa_0=\int_{-\infty}^\infty \{f(x)\}^2\,dx$. The next proposition characterises the property of infinite $*$-order.

\begin{proposition}\label{prop0}
{For $f\in \F '$ and $m\ge1$, the following statements are equivalent:
\begin{description}
\item{(i)} $\mu_{m+2}<\infty$;
\item{(ii)} $\mu_j<\infty$ for all $1\le j\le m+2$;
\item{(iii)} $\kappa_j<\infty$ and $\kappa_j=\frac{\mu_j\mu_{j+2}}{\mu_{j+1}^2}$ for all $1\le j\le m$.
\end{description}
In particular,} $f$ is of infinite $*$-order if and only if $\mu_n<\infty$, $n\ge 1$.
\end{proposition}

{\bf Proof of Proposition \ref{prop0}:\ } For each $i,n\ge 1$, {provided all terms below are finite, we have the following recursive formula

\begin{equation}\nu_{n,i}:=\int \{f^{n*}(x)\}^i\,dx=\frac{1}{\kappa_{n-1}^i}\,\nu_{n-1,i+1},\label{recursive01}\end{equation}

giving

\begin{equation}\kappa_n=\frac1{\prod_{j=0}^{n-1}\kappa_j^{n+1-j}}\,\mu_{n+2}.\label{recursive02}\end{equation}}

{(i) $\Rightarrow$ (ii) For $1\le j\le m+2$,

\begin{eqnarray*}
\mu_j&=&\int_{-\infty}^\infty \{f(x)\}^j {\bf 1}_{\{f(x)>1\}}dx+\int_{-\infty}^\infty \{f(x)\}^j {\bf 1}_{\{f(x)\le1\}}dx\\
&\le&
\int_{-\infty}^\infty \{f(x)\}^{m+2} dx+ \int_{-\infty}^\infty f(x)dx= \mu_{m+2}+1<\infty.
\end{eqnarray*}}

{(ii) $\Rightarrow$ (iii) Use \Ref{recursive02} and proceed with induction for $1\le n\le m$.}

{(iii) $\Rightarrow$ (i) By Definition~\ref{def1}, $\kappa_1<\infty$ means that $\kappa_0<\infty$. Hence (i) follows from \Ref{recursive02} with $n=m$. \qed}

\bigskip

Next we investigate the involutionary nature of the $*$-transformation.
\begin{proposition}\label{prop1}
Let $f^*$ be a \pdQ\   and assume $f^{2*}$ exists. Then
 $f^*\sim {\cal U}$ if and only if $f^{2*}\sim {\cal U}$.
\end{proposition}

{\bf Proof of Proposition \ref{prop1}:\ }{{For $r>0$, we} have

\begin{equation}\int_0^1|f^{2*}(u)-1|^r\,du=\frac1{\kappa_1^r}\int_0^1|f^*(x)-\kappa_1|^rf^*(x)\,dx.\label{xia1}
\end{equation}

If $f^*(u)\sim{\cal U}$, then $\kappa_1=1$ and (\ref{xia1}) ensures $\int_0^1|f^{2*}(u)-1|^rdu=0$, so $f^{2*}(u)\sim{\cal U}$.

Conversely, if $f^{2*}(u)\sim{\cal U}$, then using (\ref{xia1}) again gives $\int_0^1|f^*(x)-\kappa_1|^rf^*(x)\,dx=0$. Since $f^*(x)>0$ a.s.,
we have $f^*(x)= \kappa_1$ a.s. and this can only happen when $\kappa_1=1$. Thus $f^*\sim {\cal U}$, as required. \qed}.

Proposition~\ref{prop1} shows that the uniform distribution is a fixed point in the Banach space of integrable functions on [0,1] with the $L_r$ norm for any $r>0$.  It remains to show $f^{n*}$ has a limit and that the limit is the uniform distribution.
It was hoped the classical machinery for convergence in Banach spaces \cite[Ch.10]{Luen-1969} would prove useful in this regard, but the *-mapping is not a contraction. For this reason, although there are many studies of fixed point theory
in {metric and} semi-metric spaces (see, e.g., \cite{Bessenyei2017} and references therein), the fixed point Theorems \ref{xiaadd01}, \ref{coro01} and \ref{xiaadd04} shown below do not seem to be covered in these general studies. For simplicity, we use $\convlr$ to stand for the convergence in $L_r$ norm and $\convp$ for convergence in probability as $n\to\infty$.

\begin{thm}\label{xiaadd01}
 For $f\in \F '$ with infinite $*$-order, the following statements are equivalent:
\begin{description}
\item{(i)} $f^{n*}\convltwo 1$;
\item{(ii)} For all $r>0$, $f^{n*}\convlr1$;
\item{(iii)} $\frac{\mu_n\mu_{n+2}}{\mu_{n+1}^2}\to1$ as $n\to\infty$.
\end{description}
\end{thm}

\noindent{{\bf Remark} Notice that $\mu_n=\e \left\{f^*(U)^{n-1}\right\}$, $n\ge 1$, are the moments of the random variable $f^*(U)$
with $U\sim{\cal U}$, {Theorem}~\ref{xiaadd01}
says that the convergence of $\{f^{n*}:\ n\ge 1\}$ is purely determined by the moments of $f^*(U)$. This is rather puzzling because it is well known that the moments do not uniquely determine the distribution \cite[p.~227]{Feller}, meaning that different distributions with the same moments have the same converging behaviour.  {However, if $f$ is bounded, then $f^*(U)$ is a bounded random variable so its moments uniquely specify its distribution \cite[pp.~225--226]{Feller}, leading to stronger results in {Theorem}~\ref{coro01}.
}

{\bf Proof of {Theorem}~\ref{xiaadd01}:\ } It is obvious that (ii) implies (i).

(i) $\Rightarrow$ (iii): {By Proposition~\ref{prop0}, $\kappa_n=\frac{\mu_n\mu_{n+2}}{\mu_{n+1}^2}$. Now

\begin{equation}\int_0^1\{f^{n\ast}(x)-1\}^2\,dx=\kappa_n-1,\label{xiaadd02}\end{equation}

 so (iii) follows immediately.

 (iii) $\Rightarrow$ (ii): It suffices to show that $f^{n*}\convlr1$ for any integer $r\ge 4$. To this end,
 since for $a,b\ge 0$, $|a-b|^{r-2}\le a^{r-2}+b^{r-2}$, we have from (\ref{xiaadd02}) that

 \begin{equation}
 \int_0^1|f^{n*}(x)-1|^rdx\le\int_0^1(f^{n*}(x)-1)^2(f^{n*}(x)^{r-2}+1)dx=\nu_{n,r}-2\nu_{n,r-1}+\nu_{n,r-2}+\kappa_n-1,\label{xiaadd03}
 \end{equation}

 where, as before, $\nu_{n,r}=\int_0^1\{f^{n*}(x)\}^rdx$. However, {applying \Ref{recursive01} gives
 $$\nu_{n,r}=\frac{\mu_{n+r}}{\kappa_{n-1}^r\kappa_{n-2}^{r+1}\dots\kappa_0^{n+r-1}}$$
 and \Ref{recursive02} ensures
 $$\mu_{n+r}=\kappa_{n+r-2}\kappa_{n+r-3}^2\dots\kappa_0^{n+r-1},$$
 which imply}
 $$\nu_{n,r}=\kappa_{n+r-2}\kappa_{n+r-3}^2\dots \kappa_n^{r-1}\to 1$$
 as $n\to\infty$. Hence, it follows from \Ref{xiaadd03} that $ \int_0^1|f^{n*}(x)-1|^rdx\to 0$ as $n\to\infty$, completing the proof. \qed

We write $\|g\|=\sup_x|g(x)|$ for each bounded function $g$.

\begin{thm}\label{coro01}
If $f$ is bounded, then
\begin{description}
\item{(i)} {for all $n\ge 0$, $\|f^{(n+1)*}\|\le \|f^{n*}\|$} and the inequality becomes equality if and only if $f^{n*}\sim{\cal U}$;
\item{(ii)} $f^{n*}\convlr 1$ for all $r>0$.
\end{description}
\end{thm}

{\bf Proof of {Theorem}~\ref{coro01}:\ } It follows from \Ref{xiaadd02} that $\kappa_n\ge 1$ and the inequality becomes equality if and only if
 $f^{n\ast}\sim {\cal U}$.

(i) Let $Q^{n\ast}$ be the inverse of the cumulative distribution function of $f^{n\ast}$, then
$f^{(n+1)\ast}(u)=\frac{f^{n*}(Q^{n*}(u))}{\kappa_n}\le \frac{\|f^{n*}\|}{\kappa_n}$, giving $\|f^{(n+1)\ast}\|\le \frac{\|f^{n*}\|}{\kappa_n}{\le \|f^{n*}\|}$. If $f^{n\ast}\sim{\cal U}$, then
Proposition~\ref{prop1} ensures that $f^{(n+1)\ast}\sim{\cal U}$, so $\|f^{(n+1)\ast}\|=\|f^{n*}\|$. Conversely, if $\|f^{(n+1)*}\|= \|f^{n*}\|$, then $\kappa_n=1$, so $f^{n\ast}\sim {\cal U}$.

(ii) It remains to show that $\kappa_n\to 1$ as $n\to\infty$. In fact, if $\kappa_n\not\to 1$, since $\kappa_n\ge 1$, there exist a $\delta>0$ and a subsequence $\{n_k\}$ such that $\kappa_{n_k}\ge 1+\delta$, which implies

\begin{equation}\frac{\mu_{n_k+2}}{\mu_{n_k+1}}=\prod_{i=0}^{n_k}\kappa_i\ge (1+\delta)^k\to\infty\mbox{ as }k\to\infty.\label{xia3}\end{equation}

However, $\frac{\mu_{n_k+2}}{\mu_{n_k+1}}\le \|f\|<\infty$, which contradicts (\ref{xia3}).  \qed

\begin{thm}\label{xiaadd04}
For $f\in \F '$ with infinite $*$-order such that $\{\mu_n\mu_{n+2}\mu_{n+1}^{-2}:\ n\ge 1\}$ is a bounded sequence, then the following statements are equivalent:
\begin{description}
\item{(i*)} $f^{n*}\convp1$;
\item{(ii)} For all $r>0$, $f^{n*}\convlr1$;
\item{(iii)} $\mu_n\mu_{n+2}\mu_{n+1}^{-2}\to1$ as $n\to\infty$.
\end{description}
\end{thm}

{\bf Proof of {Theorem}~\ref{xiaadd04}:\ } It suffices to show that (i*) implies (iii). Recall that $\kappa_n=\mu_n\mu_{n+2}\mu_{n+1}^{-2}$, for each subsequence
$\{\kappa_{n_k}\}$, there exists a converging sub-subsequence $\{\kappa_{n_{k_i}}\}$ such that $\kappa_{n_{k_i}}\to b$ as $i\to\infty$. It remains to show that $b=1$. To this end, for $\delta>1$, we have

\begin{eqnarray}
&&\int_0^1\left|f^{(n_{k_i}+1)*}(x)-1\right|{\bf 1}_{\left\{\left|f^{(n_{k_i}+1)*}(x)-1\right|\le \delta\right\}}dx\nonumber\\
&&=\frac1{\kappa_{n_{k_i}}}\int_0^1\left|f^{(n_{k_i})*}(x)-\kappa_{n_{k_i}}\right|f^{(n_{k_i})*}(x){\bf 1}_{\left\{\left|f^{n_{k_i}*}(x)-\kappa_{n_{k_i}}\right|\le \delta\kappa_{n_{k_i}}\right\}}dx.\label{xiaadd05}
\end{eqnarray}

(i*) ensures that  $$\left|f^{(n_{k_i}+1)*}-1\right|\convp0, \ f^{n_{k_i}*}\left|f^{n_{k_i}*}-\kappa_{n_{k_i}}\right|\convp |1-b|, \
{\bf 1}_{\left\{\left|f^{n_{k_i}*}(x)-\kappa_{n_{k_i}}\right|\le \delta\kappa_{n_{k_i}}\right\}}\convp 1$$ as $i\to\infty$, so applying the bounded convergence theorem to both sides of \Ref{xiaadd05} to get $0=|1/b-1|$, i.e., $b=1$. \qed

\subsection{Examples of convergence to uniformity}\label{sec:examplesconv}

The main results in section~\ref{sec:convunif} cover all the standard distributions with infinite $*$-order in
\cite{J-K-B-1994}, \cite{J-K-B-1995}. In fact, as observed in the Remark after {Theorem}~\ref{xiaadd01} that the convergence to uniformity is purely determined by the moments of $f^*(U)$ with $U\sim{\cal U}$,
we have failed to construct a density such that $\{f^{n*}:\ n\ge 1\}$ does not converge to the uniform distribution. Here we give a few examples to show that the main results in section~\ref{sec:convunif} are indeed very convenient to use.

\subsubsection*{Example 1: Power function family.}
From Table~\ref{table1}  the Power$(b)$ family  has density $f_b(x)= bx^{b-1},~0<x<1$, so it is of infinite $*$-order if and only if $b\ge 1$. As $f_b$ is
bounded for $b\ge 1$, {Theorem}~\ref{coro01}
ensures that $f^{n*}_b$ converges to the uniform in $L_r$ for any $r>0$.

\subsubsection*{Example 2: Exponential distribution.}

Suppose $f(x)=e^x,\; ~x<0$.
{Then $f^*(u)=2u$, $0<u<1${,} which belongs to  the Power(2) distribution;
and so by Example~1,} $f$ is bounded, so {Theorem}~\ref{coro01} says that $f^{n*}$ converges to the uniform distribution as $n\to \infty.$ By symmetry,
the same result holds for $f(x)=e^{-x},\; ~x>0$.

\subsubsection*{Example 3: Pareto distribution.}

The Pareto($a$) family, with $a>0${,} has $f_a(x)= ax^{-a-1}$ for $ x> 1 $, which is bounded, so an application of {Theorem}~\ref{coro01} yields that the sequence $\{f^{n*}_a\}_{n\geq 1}$ converges to the uniform distribution as $n\to \infty.$

\subsubsection*{Example 4: Cauchy distribution.}

The \pdQ\ of  the Cauchy density is given by $f^*(u)=2\sin ^2(\pi u)$, $0<u<1$, see Table~\ref{table1}; it retains the bell shape of $f$.  It follows that
 $F^*(t)=t-\sin (2\pi t)/(2\pi ),$ for $0<t<1$. It seems impossible to obtain an analytical form of $f^{n*}$ for $n\ge 2$. However, as $f$ is bounded, using {Theorem}~\ref{coro01}, we can conclude that $f^{n*}$ converges to the uniform distribution as $n\to\infty$.

\subsubsection*{Example 5: Normal distribution.}

Although it is possible to obtain $\{f^{n*}\}$ by induction and then derive directly that $f^{n*}$ converges to the uniform distribution as $n\to\infty$, one can easily see that the pdf is bounded and so {Theorem}~\ref{coro01} can be employed to get the same conclusion.

\subsubsection*{Example 6:} Let $f(x)=-\ln x$, $x\in(0,1)$, then $\mu_n=n!$ and $\kappa_n=\frac{n+2}{n+1}\to1$ as $n\to\infty$, so we have from {Theorem}~\ref{xiaadd01} that for any $r>0$, $f^{n*}$ converges in $L_r$ norm to constant 1 as $n\to\infty$.



\section{Testing for uniformity}\label{sec:testing}

The larger the value of $d({\cal U},f^*)$, the easier it should be to discriminate between $\cal U$ and $f^*$.
This is indeed the case, as we now demonstrate.

\subsection{Power of tests for detecting non-uniform shapes}\label{sec:power}

 The connection between Kullback-Leibler divergences and the Neyman-Pearson Lemma is well-known, see text and references in \cite{ecopas-2006}. The following material on power functions for tests for uniformity, while contextual to comparing shapes, may prove useful  in other situations.

Given a random sample of $m$ independent, identically distributed (\iid) variables, each from a  distribution with density $f$, it is feasible to carry out a nonparametric test of uniformity by estimating the \pdQ \ with a kernel density estimator $\widehat f_m^*$ and comparing it with the uniform density on [0,1] using any one of a number of  metrics. Consistent estimators $\widehat f_m^*$ for $f^*$ based on normalized reciprocals of the quantile density  estimators derived in \cite{PRST-2016a} are available and described in \citet[Section~2]{S-2017}. However an investigation into such omnibus nonparametric testing
 procedures, and comparison with other kernel density based techniques {in} \cite{bow-1992,fan-1994} {and} \cite{pav-2015}, is beyond the scope of this work.

 Given the large number of tests for uniformity that are available, see text and references in \cite{step-2006},
 one may well ask why introduce new ones? The {\em doyen} of goodness-of-fit testing \cite{step-2006} provides the answer:

 \begin{quote}\em
  Since transformations are often used to produce
a set of uniforms, it might be appropriate
to conclude with some cautionary remarks
on when uniformity is not to be expected.
This will be so, for example, when the U
set is derived from the PIT  and when some
parameters, unknown in the distribution, are
replaced by estimates. In this situation, even
when the estimates are efficient, the U set
will be superuniform, giving much lower values
of, say, EDF statistics, than if the set were
uniform; this remains so even as the sample
size grows bigger.
 \end{quote}

In practice this means that if a test for uniformity is preceded by a probability integral transformation (PIT) including parameter estimates, then the actual levels of such tests will not be those nominated unless
(often complicated and model specific) adjustments are made. Examples are in \cite{lock-1986} and \cite{schsch-1997}.

 In this section we study the simpler problem of
 testing the null hypothesis of uniformity $H_0: f^*(u)=1$ for all $0<u<1$ (denoted {$\cal U$}) against  a specified alternative $H_1: f^*=f^*_1$. This test will  give us a standard by which to judge the power of
 any future nonparametric test when the specific alternative holds.
Given a vector of $\bX =(X_1,\dots ,X_m)$ of \iid \;variables from $f^*$, and realized values $\bx =(x_1,\dots ,x_m)$ the  Neyman-Pearson (NP) test rejects $H_0$ in favor of $H_1$ when the product $\prod_{1}^{m} f_1^*(x_i)$ is large, or equivalently, when

\begin{equation}\label{eqn:NPtestexact}
    l_\bx = \sum_{1}^{m}\ln(f_1^*(x_i)) \geq c_{m,\alpha}~,
\end{equation}

where $c_{m,\alpha}$ is chosen to achieve level $\alpha $. In general this critical point will be difficult to
determine, but for the normal \pdQ \  alternative $f_1^*(u)= 2\sqrt \pi\,\varphi(z_u) $ we have:
$l_\bx = m\ln (\sqrt 2)-\frac{1}{2}\,\sum_{1}^{m}\{\Phi ^{-1}(x_i)\}^2$, where $\varphi$ is the density of the standard normal.
Under the null hypothesis $\sum_{1}^{m}\{\Phi ^{-1}(U_i)\}^2\sim\chi ^2_m$, so an equivalent test
to (\ref{eqn:NPtestexact}) would reject the null hypothesis of uniformity in favor of normality
when $\sum_{1}^{m}\{\Phi ^{-1}(x_i)\}^2\leq \chi ^2_m(\alpha ).$ This is the most powerful level-$\alpha $ test of these simple hypotheses based on $m$ observations.

Returning to a general alternative \pdQ\  we can find asymptotically most powerful level-$\alpha $ tests
based on the fact that $l_\bX $ is a sum of \iid \;random variables with common mean $\mu _0 = \e _0[\ln(f_1^*(U))] $ and variance $\sigma^2 _0 = \var _0[\ln(f_1^*(U))] $, which we assume are finite. By virtue of a CLT, the
 large-sample NP test rejects $H_0$ at asymptotic level $\alpha $ when

\begin{equation}\label{eqn:NPtest}
    { (l_\bx/\sqrt{m} -\sqrt{m}\mu _0 )/\sigma _0} \geq z_{1-\alpha }~.
\end{equation}

To obtain an expression for the large sample power of this test,
let $\mu _1 = \e _1[\ln(f_1^*(X))] $ and  $\sigma^2 _1 = \var _1[\ln(f_1^*(X))] $, again assumed to be finite; then the asymptotic power of
the test (\ref{eqn:NPtest}) against $f^*$ based on $m$ observations is readily found to be:

\begin{equation}\label{eqn:NPpower}
\Pi _m (f_1^*)=\Phi \left(\sqrt m \;\frac{(\mu _1-\mu _0)}{\sigma _1}+z_\alpha \;\frac {\sigma _0}{\sigma _1}\right )~.\end{equation}

We need $\alpha $, $m$, $l_\bx $,  $\mu_0$, $\sigma^2 _{0}$ to carry out the  test (\ref{eqn:NPtest}); and we
also need $\mu_1$, $\sigma^2 _1$ to compute the asymptotic power (\ref{eqn:NPpower}).
Notice that the distances from the origin in Figure~\ref{fig1} of  Section~\ref{sec:defns} are based on the  directed divergences $I({\cal U}:f_1^*)=-\e _0[\ln (f_1^*(X))]=-\mu _0$ and $I(f_1^*:{\cal U})=\e _1[\ln (f_1^*(X))]=\mu _1$, so the symmetrized divergence  KLD is $J({\cal U},f_1^*)= \mu _1-\mu _0.$ Thus the power function (\ref{eqn:NPpower}) is non-decreasing in the KLD or its square root, the distance $d({\cal U},f^*_1)$ between null and alternative.

\begin{table}[t!]
\begin{small}
\begin{center}
\caption{{\bf Quantiles of some distributions, their \pdQs\ and quantities relevant to the asymptotic power
function (\ref{eqn:NPpower}).} \em In general, we denote $x_u=Q(u)=F^{-1}(u)$, but for the normal $F=\Phi $ with density $\varphi $, we use $z_u=\Phi ^{-1}(u)$.
The logistic quantile function is only defined {for} $u\leq 0.5$ but it is symmetric about $u=0.5.$  LN represents the
standard lognormal distribution. The quantile function
for the Pareto is for the Type~II distribution with shape $a=1$, and the \pdQ \  is the same for Type~I and Type~II Pareto models.  \label{table1}}
\vspace{.1cm}
\setlength{\tabcolsep}{1mm}
\renewcommand{\arraystretch}{1.2}
\begin{tabular}{lllccccccc}
\hline
           &              &                  & $I({\cal U}:f^*)$ & $I(f^*:{\cal U})$ &&& $J({\cal U},f^*)$ && \\
           & $\quad Q(u)$ &  $\quad f^*(u)$  & $-\mu _0$ & $\mu _1$ &  $\sigma _0$ & $\sigma _1 $ & $\mu _1-\mu _0$ &
           $\frac {\mu _1-\mu _0}{\sigma _1} $ &  $\frac{\sigma _0}{\sigma _1}$ \\
\cline{2-10}\\[-.3cm]
Normal   & $z_u$          &  $2\sqrt \pi\,\varphi(z_u)$     & 0.153 & 0.097 & 0.707 & 0.354 & 0.250 & 0.707 & 2.000  \\
Logistic & $\ln(u/(1-u))$ &  $6u(1-u)$                      & 0.208 & 0.125 & 0.843 & 0.393 & 0.333 & 0.848 & 2.143  \\
Laplace  & $\ln (2u),~u\leq 0.5$&$2\,\min\{u,1-u\}$         & 0.307 & 0.193 & 1.000 & 0.500 & 0.500 & 1.000 & 2.000  \\
$t_2$    & $\frac{2u-1}{\{2u(1-u)\}^{1/2}}$ &$\frac{2^7\{u(1-u)\}^{3/2}}{3\pi }$ & 0.391 & 0.200 & 1.264 & 0.463 &
0.591 & 1.276 & 2.728 \\
Cauchy   & $\tan\{\pi(u-0.5)\}$ & $2\sin ^2(\pi u)$         & 0.693 & 0.307 & 1.814 & 0.538 & 1.000 & 1.857 & 3.369 \\
Exp.     & $-\ln (1-u)$   &  $2(1-u)$                       & 0.307 & 0.193 & 1.000 & 0.500 & 0.500 & 1.000 & 2.000 \\
Gumbel   & $-\ln(-\ln (u))$ &  $-4u\ln(u)$                  & 0.191 & 0.116 & 0.803 & 0.381 & 0.307 & 0.806 & 2.109 \\
LN     & $e^{z_u}$ & $\frac {2\sqrt \pi\;}{e^{1/4}}\,\varphi (z_u)\,e^{-z_u}$ & 0.403 & 0.222 & 1.225 & 0.500 & 0.625 & 1.250 &  2.449 \\
Pareto &  $(1-u)^{-1}$ & $3\,(1-u)^{2}$       & 0.901 & 0.432 & 2.000 & 0.667 & 1.333 & 2.000 & 3.000 \\
\hline
 \end{tabular}
\end{center}
\end{small}
\end{table}

Some examples of $\mu _0$, $\mu _1$, $\mu _1-\mu _0$, $(\mu _1-\mu _0)/\sigma _1$ and $\sigma _0/\sigma _1$
 are given in Table~\ref{table1}. Note the particularly
simple values for the symmetrized divergence $J({\cal U},f_1^*)=\mu _1-\mu _0.$ Distributions with shapes \lq visually far\rq \  from uniformity have large values in the two right-most  columns, so that the test will more easily detect them. For the normal alternative, the asymptotic power (\ref{eqn:NPpower}) at level $\alpha $ is $\Pi _m (f_1^*)=\Phi (\sqrt {m/2} +2\,z_\alpha ),$ which exceeds $\alpha $ when  $m> 2z_{1-\alpha }^2.$

There are other situations where the asymptotic power functions are monotone increasing functions of the KLD
with many one-sample examples in \cite{KMS-2008,KMS-2010,KMS-2014}, one-parameter
 exponential families \cite{M-S-2012}, two-sample binomial tests \cite{prst-2014} and
non-central chi-squared and non-central F families arising {in} tests for equivalence \cite{M-S-2016}. This {reveals a} general phenomenon but no meta-theorem containing these results is yet available.

\begin{figure}[t!]
\begin{center}
\includegraphics[scale=.8]{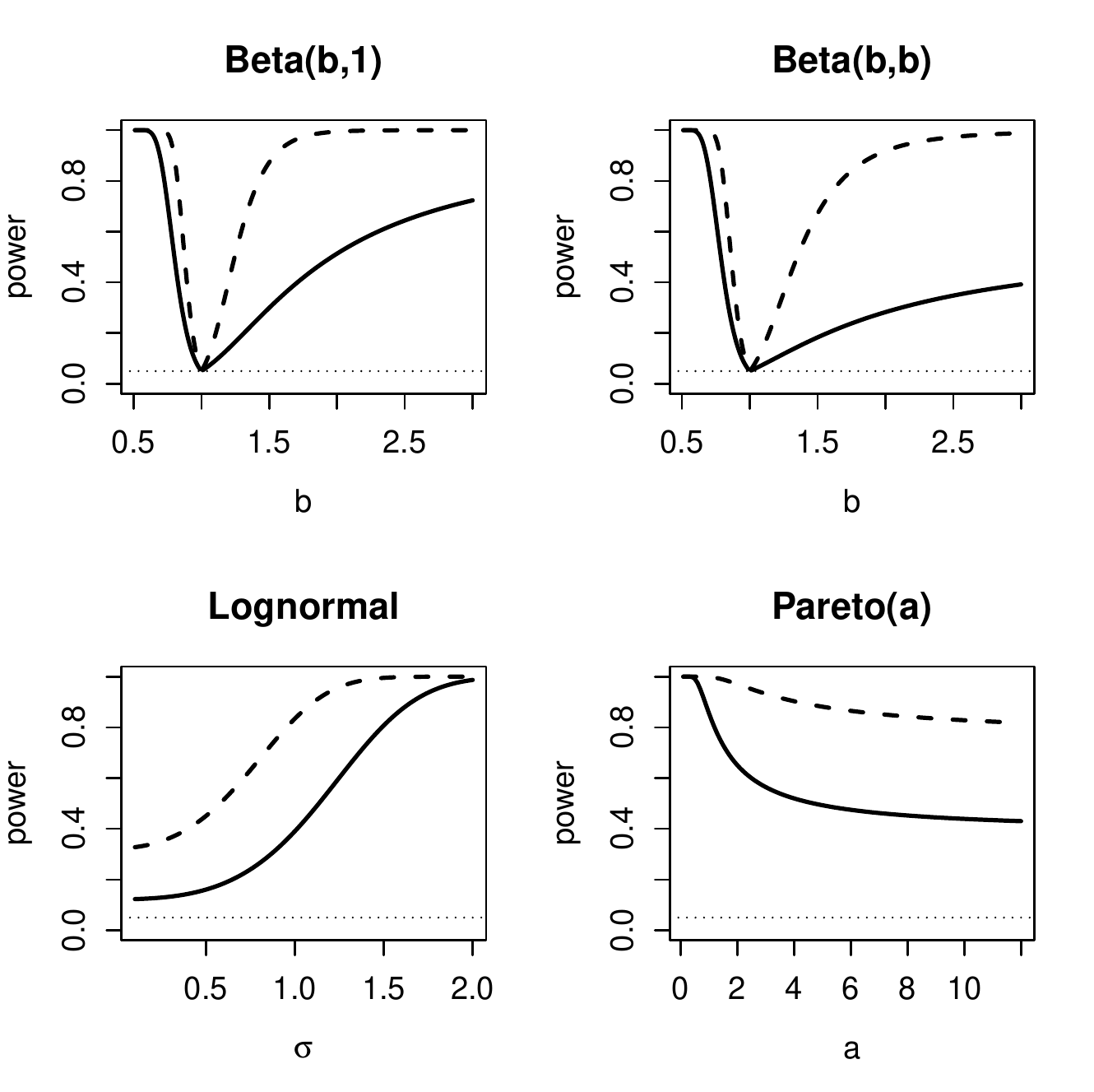}
\caption{\small The top left plot shows the asymptotic power function for level 0.05 tests of uniformity against alternative the Beta(b,1) \pdQ\ when $m=25$ (solid line) and $m=100$ (dashed line). The plots for the {symmetric}
Beta distribution are on its right, with the same sample   sizes; these power functions show that it is harder to
detect the symmetric ones for all $b\neq 1.$   In the bottom two plots, the sample sizes are much smaller,
$m=9 $ solid lines and $m=16$ dashed lines.
  \label{fig2}}
 \end{center}
 \end{figure}

\subsection{Examples of power functions for shape families}\label{sec:shapeexs}
 The power functions of testing uniformity against the \pdQs \  of  four shape families,
 are shown in {Figure~\ref{fig2}}. The first two models, power function model Beta($b,1$) and the symmetric  Beta($b,b$)
 model for $b> 0.5$ contain the null hypothesis. Their respective  power functions (\ref{eqn:NPpower}) for a level 0.05 test of uniformity based on $m=25$ and 100 observations are shown in the top two plots. Similar plots for alternatives Lognormal($\sigma$) and Pareto($a$) families  are also shown for much smaller sample sizes 9 and 16
 indicating that small samples will likely detect these alternative shapes.

The plots in {Figure~\ref{fig2}} require the null and alternative means and variances of the
test statistic, and were obtained by numerical integration. In the case of the Beta($b$,1)  model exact
results are derived as follows. the quantile function is $Q(b)=u^{1/b}$ and for  $b>1/2$ its density
is square integrable, leading to the \pdQ \  $f_1^*(u)=(2-\frac {1}{b})\,u^{1-\frac {1}{b}}.$  The log-likelihood
for one observation $X=x$ required in (\ref{eqn:NPtest}) is $l_x(b)=\ln(2-1/b)+(1-1/b)\ln(x)$. Thus
 $\mu _0(b) =\e _0[l_\bX (b)]=\ln(2-1/b)+1/b-1$  and  $\mu _1(b)-\mu _0(b) =(1-1/b)^2/(2-1/b).$  Further
$\sigma _0(b) = \{\var _0[l_\bX (b)]\}^{1/2}=|1-1/b| $ and $\sigma _1(b) = \{\var _1[l_\bX (b)]\}^{1/2}=|1-1/b|/(2-1/b).$
Hence $(\mu _1(b)-\mu _0(b))/\sigma _1(b) =|1-1/b|$ and $\sigma _0(b)/\sigma _1(b) = 2-1/b.$
The asymptotic power function (\ref{eqn:NPpower}) is therefore
$\Pi _m (f_1^*(b))=\Phi (\sqrt {m}\; |1-1/b| +(2-1/b)\,z_\alpha )$ for $ b>1/2.$

\section{Summary and Discussion}\label{sec:summary}

The \pdQ\, transformation from a density function $f$ to $f^*$ extracts the important information of $f$ such as its asymmetry and tail behaviour and ignores
the less critical information such as gaps, location and scale and thus provides a powerful tool in studying the shapes of density functions.

We found the directed divergences from uniformity of the \pdQs\ of many standard location-scale families
and used them to make a map locating each shape family relative to others and giving its distance from uniformity. We also found the most powerful tests of uniformity against alternative shapes and showed that their power functions are monotone increasing in the distances from the origin on the map.

In terms of the limiting behaviour of repeated applications of the \pdQ\ mapping, when the density function $f$ is bounded, we showed that each application lowers its modal height and hence the resulting
density function $f^*$ is closer to the uniform density than $f$. Furthermore, we established a necessary and sufficient condition for
$f^{n*}$ converging in $L_2$ norm to the uniform density, giving a positive answer to a conjecture raised in
\cite{S-2017}. In particular, if $f$ is bounded, we proved that $f^{n*}$ converges in $L_r$ norm to the uniform density for any $r>0$.
The fixed point {theorems} can be interpreted as follows. As we repeatedly apply the \pdQ\ transformation, we keep losing information about the  shape of the original $f$ and will eventually exhaust the information,
leaving nothing in the limit, as represented by the uniform density, which means no points carry more information than other points. Thus the \pdQ\ transformation plays
a similar role to the difference operator in time series analysis where repeated applications of the difference operator to a time series with polynomial component  lead to a white noise with a constant power spectral density \cite[p.~19]{Brockwell-2009}.

We conjecture that every almost surely positive density $g$ on $[0,1]$ is a \pdQ\ of a density function, hence uniquely represents a location-scale family. This is equivalent to saying that  there exists a density function $f$ such that $g=f^*$. When $g$ satisfies $\int_0^1\frac1{g(t)}dt<\infty$, one can show that
the \cdf $F$ of $f$ can be uniquely (up to location-scale parameters) represented as $F(x)=H^{-1}(H(1)x)$, where $H(x)=\int_0^x\frac1{g(t)}dt$ (Professor A.D. Barbour, personal communication).
The condition $\int_0^1\frac1{g(t)}dt<\infty$ is equivalent to saying that $f$ has bounded support and it is certainly not necessary, e.g., $g(x)=2x$ for $x\in [0,1]$ and $f(x)=e^{x}$ for $x<0$ (see Example~2 in Section~\ref{sec:examplesconv}).

In summary, the study of shapes of probability densities is facilitated by composing them with their own
quantile functions, which puts them on the same finite support where they are absolutely continuous with respect to
{Lebesgue} measure, and thus amenable to metric and semi-metric comparisons. In addition, we showed that further applications
of this transformation, which intuitively reduces information and increases the relative entropy, is generally valid
but requires a non-standard approach for proof.  Similar results are likely to be obtainable in
the multivariate case.   Further research could investigate the relationship between relative entropy and tail-weight or distance from the class of symmetric \pdQs.

{\em Acknowledgments: \ The authors thank Professor P.J. Brockwell for helpful commentary on an earlier version of this manuscript. The research by Professor Aihua Xia is supported by ARC Discovery Grant DP150101459}


\begin{thebibliography}{}

\bibitem[\protect\citename{Bessenyei \& P\'{a}les, }2017]{Bessenyei2017}
{\sc Bessenyei, M., \& P\'{a}les, Z.} 2017.
\newblock A contraction principle in semimetric spaces.
\newblock {\em J. nonlinear convex anal.}, {\bf 18}(3), 515--524.

\bibitem[\protect\citename{Bowman, }1992]{bow-1992}
{\sc Bowman, A.W.} 1992.
\newblock Density based tests for goodness-of-fit.
\newblock {\em J. {S}tatist. {C}omp. \& {S}im.}, {\bf 40}, 1--13.

\bibitem[\protect\citename{Brockwell \& Davis, }2009]{Brockwell-2009}
{\sc Brockwell, P.J., \& Davis, R.A.} 2009.
\newblock {\em Time {S}eries: {T}heory and {M}ethods}.
\newblock Springer-Verlag.

\bibitem[\protect\citename{Eguchia \& Copas, }2006]{ecopas-2006}
{\sc Eguchia, S., \& Copas, J.} 2006.
\newblock Interpreting kullback–leibler divergence with the neyman–pearson
  lemma.
\newblock {\em J. {M}ultiv. {A}nal.}, {\bf 97}, 2034--2040.

\bibitem[\protect\citename{Fan, }1994]{fan-1994}
{\sc Fan, Y.} 1994.
\newblock Testing the goodness of fit of a parametric density function by
  kernel method.
\newblock {\em Econometric {T}heory}, {\bf 10}, 316--356.

\bibitem[\protect\citename{Feller, }1971]{Feller}
{\sc Feller, W.} 1971.
\newblock {\em An introduction to probability theory and its applications,
  vol.~2}.
\newblock New York: John Wiley \& Sons, Inc.

\bibitem[\protect\citename{Freimer {\em et~al.}, }1988]{fmkl-1988}
{\sc Freimer, M., Mudholkar, G.S., Kollia, G., \& Lin, C.T.} 1988.
\newblock A study of the generalized {T}ukey lambda family.
\newblock {\em Commun. {S}tatist. - a}, {\bf 17}, 3547--3567.

\bibitem[\protect\citename{Johnson {\em et~al.}, }1994]{J-K-B-1994}
{\sc Johnson, N.L., Kotz, S., \& Balakrishnan, N.} 1994.
\newblock {\em Continuous univariate distributions}.
\newblock  Vol. 1.
\newblock New York: John Wiley \& Sons.

\bibitem[\protect\citename{Johnson {\em et~al.}, }1995]{J-K-B-1995}
{\sc Johnson, N.L., Kotz, S., \& Balakrishnan, N.} 1995.
\newblock {\em Continuous univariate distributions}.
\newblock  Vol. 2.
\newblock New York: John Wiley \& Sons.
\newblock ISBN 0-471-58494-0.

\bibitem[\protect\citename{Kulinskaya {\em et~al.}, }2008]{KMS-2008}
{\sc Kulinskaya, E., Morgenthaler, S., \& Staudte, R.G.} 2008.
\newblock {\em Meta {A}nalysis: a {G}uide to {C}alibrating and {C}ombining
  {S}tatistical {E}vidence}.
\newblock Wiley {S}eries in {P}robability and {S}tatistics.
\newblock Chichester: John {W}iley \& {S}ons.

\bibitem[\protect\citename{Kulinskaya {\em et~al.}, }2010]{KMS-2010}
{\sc Kulinskaya, E., Morgenthaler, S., \& Staudte, R.G.} 2010.
\newblock Variance stabilizing the difference of two binomial proportions.
\newblock {\em {A}mer. {S}tatist.}, {\bf 64}(4), 350--356.

\bibitem[\protect\citename{Kulinskaya {\em et~al.}, }2014]{KMS-2014}
{\sc Kulinskaya, E., Morgenthaler, S., \& Staudte, R.G.} 2014.
\newblock Combining {S}tatistical {E}vidence.
\newblock {\em Int. {S}tatist. {I}nst. {R}ev.}, {\bf 82}(2), 214--242.

\bibitem[\protect\citename{Kullback, }1968]{Kullback-1968}
{\sc Kullback, S.} 1968.
\newblock {\em Information {T}heory and {S}tatistics}.
\newblock Mineola, NY: Dover.

\bibitem[\protect\citename{Kullback \& Leibler, }1951]{kl-1951}
{\sc Kullback, S., \& Leibler, R.A.} 1951.
\newblock On information and sufficiency.
\newblock {\em Ann. {M}ath. {S}tatist.}, {\bf 22}, 79--86.

\bibitem[\protect\citename{Lockhart {\em et~al.}, }1986]{lock-1986}
{\sc Lockhart, R.A., O'Reilly, F.J., \& Stephens, M.A.} 1986.
\newblock Tests of fit based on normalized spacings.
\newblock {\em J. {R}. {S}tatist. {S}oc. {B}}, {\bf 48}, 344--352.

\bibitem[\protect\citename{Luenberger, }1969]{Luen-1969}
{\sc Luenberger, D.G.} 1969.
\newblock {\em {Optimization by Vector Space Methods}}.
\newblock New York, NY: Wiley.

\bibitem[\protect\citename{Morgenthaler \& Staudte, }2012]{M-S-2012}
{\sc Morgenthaler, S., \& Staudte, R.G.} 2012.
\newblock Advantages of variance stabilization.
\newblock {\em Scand. {J}. {S}tatist.}, {\bf 39}, 714--728.

\bibitem[\protect\citename{Morgenthaler \& Staudte, }2016]{M-S-2016}
{\sc Morgenthaler, S., \& Staudte, R.G.} 2016.
\newblock Indicators of evidence for equivalence.
\newblock {\em Entropy}, {\bf 18}, 291.

\bibitem[\protect\citename{Parzen, }1979]{par-1979}
{\sc Parzen, E.} 1979.
\newblock Nonparametric statistical data modeling.
\newblock {\em J. {A}mer. {S}tatist. {A}ssoc.}, {\bf 7}, 105--131.

\bibitem[\protect\citename{Pavia, }2015]{pav-2015}
{\sc Pavia, J.M.} 2015.
\newblock Testing goodness-of-fit with the kernel density estimator: Gofkernel.
\newblock {\em J. {S}tatist. {S}oft.}, {\bf 66}, 1--27.

\bibitem[\protect\citename{Prendergast \& Staudte, }2014]{prst-2014}
{\sc Prendergast, L.A., \& Staudte, R.G.} 2014.
\newblock Better than you think: interval estimators of the difference of
  binomial proportions.
\newblock {\em J. {S}tatist. {P}lann. {I}nference}, {\bf 148}, 38--48.

\bibitem[\protect\citename{Prendergast \& Staudte, }2016]{PRST-2016a}
{\sc Prendergast, L.A., \& Staudte, R.G.} 2016.
\newblock Exploiting the quantile optimality ratio in finding confidence
  intervals for a quantile.
\newblock {\em {Stat}}, {\bf 5}(1), 70--81.

\bibitem[\protect\citename{Schader \& Schmid, }1997]{schsch-1997}
{\sc Schader, M., \& Schmid, F.} 1997.
\newblock Power of tests for uniformity when limits are unknown.
\newblock {\em J. {A}ppl. {S}tatist.}, {\bf 24}, 193--205.

\bibitem[\protect\citename{Staudte, }2017]{S-2017}
{\sc Staudte, R.G.} 2017.
\newblock The shapes of things to come: probability density quantiles.
\newblock {\em Statistics}, {\bf 51}, 782--800.

\bibitem[\protect\citename{Stephens, }2006]{step-2006}
{\sc Stephens, M.S.} 2006.
\newblock Uniformity, tests of.
\newblock {\em Pages  1--8 of:} {\em Encyclopedia of {S}tatistical {S}ciences},
   vol. 53.
\newblock John {W}iley \& {S}ons.
\newblock DOI: 10.1002/0471667196.ess2810.pub2.

\bibitem[\protect\citename{Team, }2008]{R}
{\sc Team, R Development~Core}. 2008.
\newblock {\em R: A language and environment for statistical computing}.
\newblock R Foundation for Statistical Computing, Vienna, Austria.
\newblock {ISBN} 3-900051-07-0.

\end{thebibliography}

\end{document}